\documentclass{amsart}
\usepackage{amsmath}
\usepackage{amsfonts}
\usepackage{amssymb}
\usepackage{amscd}
\usepackage{verbatim}
\usepackage[arrow,matrix,cmtip]{xy}

\def\cprime{$'$} 
\newcommand{\bpr}{\begin{proof}}
\newcommand{\epr}{\end{proof}}

\newcommand{\lm}{\lim_{n \ra \infty}}

\newcommand{\ext} {\operatorname{Ext}}

\newcommand{\lra} {\longrightarrow}

\newcommand{\ra}{\rightarrow}

\newcommand{\mc}{\mathcal}
\newcommand{\mf}{\mathfrak}
\newcommand{\cha}{\operatorname{char}}
\newcommand{\mb}{\mathbb}

\newcommand{\pgl}{\operatorname{PGL}}

\newcommand{\uext}{\underline{\operatorname{Ext}}}

\newcommand{\uhom}{\underline{\operatorname{Hom}}}
\renewcommand{\hom}{\operatorname{Hom}}

\newcommand{\wt}{\widetilde}

\newcommand{\coH}{\operatorname{H}}
\newcommand{\ucoH}{\underline{\operatorname{H}}}

\newcommand{\tor}{\operatorname{Tor}}
\newcommand{\utor}{\underline{\operatorname{Tor}}}

\newcommand{\cd}{\operatorname{cd}}

\newcommand{\Proj}{\operatorname{-Proj}}

\newcommand{\cproj}{\operatorname{proj}}

\newcommand{\rProj}{\operatorname{Proj-}}

\newcommand{\gd}{\operatorname{gd}}
\newcommand{\Gr}{\operatorname{-Gr}}
\newcommand{\gr}{\operatorname{-gr}}
\newcommand{\rQch}{\operatorname{Qch}}
\newcommand{\Qgr}{\operatorname{-Qgr}}

\newcommand{\rQgr}{\operatorname{Qgr-}}

\newcommand{\Tors}{\operatorname{-Tors}}

\newcommand{\aut}{\operatorname{Aut}}

\newcommand{\segre}{\overset{s}{\otimes}}

\newtheorem{theorem}{Theorem}[section]
\newtheorem{lemma}[theorem]{Lemma}

\newtheorem{proposition}[theorem]{Proposition}

\theoremstyle{definition}
\newtheorem{hypothesis}[theorem]{Hypothesis}

\newtheorem{remark}[theorem]{Remark}

\numberwithin{equation}{section}
\makeatletter                      
\let\c@equation\c@theorem          
\makeatother     

\begin{document}
\title{Idealizer Rings and Noncommutative Projective Geometry}
\author{Daniel Rogalski}
\address{Department of Mathematics, MIT, Cambridge, MA 02139-4307. } 
\email{rogalski@math.mit.edu}
\keywords{noetherian graded ring, noncommutative projective geometry, idealizer ring}
\subjclass{Primary 16W50, Secondary 14A22}
\thanks{The author was supported in part by NSF grants NSF-DMS-9801148 and NSF-DMS-0202479, and 
a Clay Mathematics Institute Liftoff Fellowship. \hfill}

\begin{abstract}
We study noetherian graded idealizer rings which have very different behavior on the 
right and left sides. In particular, we construct noetherian graded algebras $T$ over an algebraically 
closed field $k$ with the following properties:  $T$ is left but not right strongly 
noetherian; $T \otimes_k T$ is left but not right noetherian and $T \otimes_k T^{op}$ 
is noetherian; the left noncommutative projective scheme $T\Proj$ is different from 
the right noncommutative projective scheme $\rProj T$; and 
$T$ satisfies left $\chi_d$ for some $d \geq 2$ yet fails right $\chi_1$.
\end{abstract}
\maketitle

\section{Introduction}
As a general principle, rings 
which are both left and right noetherian are expected to have rather 
symmetric properties on their left and the right sides.  The 
theme of this paper is to show that  
such intuition fails quite utterly
for certain properties which are important in the theory of noncommutative
projective geometry.  
Our main result is the following.
\begin{theorem}(Theorem~\ref{main theorem redux})
\label{main theorem} For any integer $d \geq 2$, there exists a connected finitely 
presented graded noetherian $k$-algebra $T$, where $k$ is an algebraically closed 
field, such that 
\begin{enumerate}
\item $T$ is strongly left noetherian, but not strongly right noetherian;
\item $T \otimes_k T$ is left but not right noetherian, while $T \otimes_k T^{op}$ is 
noetherian;
\item the noncommutative projective schemes
$T\Proj$ and $\rProj T$ have equivalent underlying categories, but non-isomorphic 
distinguished objects; and  
\item $T$ satisfies $\chi_{d-1}$ but not $\chi_d$ on the left, yet $T$ fails $\chi_1$ on the 
right.
\end{enumerate} 
\end{theorem}
\noindent In the remainder of the introduction, we 
will define and briefly discuss all of the relevant terms in the statement of the theorem 
and indicate how the ring $T$ is constructed.  For a more detailed 
introduction to the theory of noncommutative geometry which motivates the 
study of these properties, see the survey article \cite{SVdB}.

If $R$ is a $k$-algebra, then 
$R$ is called \emph{strongly left (right) noetherian} if $R \otimes_k B$ is left 
(right) noetherian for every commutative noetherian $k$-algebra $B$.  
The study of the strong noetherian condition for graded rings in particular has recently become 
important because of the appearance of this property in the hypotheses of 
several theorems in noncommutative geometry.  Most notably, Artin and Zhang showed 
that if $A$ is a strongly noetherian graded $k$-algebra, then 
the set of graded $A$-modules with a given Hilbert function 
is parametrized by a projective scheme \cite{AZ}.  
It is not a priori obvious that any noetherian finitely generated $k$-algebra which is not
strongly noetherian should exist; in \cite{RS}, Resco and Small gave the first (ungraded) such example.  
More recently, the author showed 
that there exist noncommutative noetherian graded rings which are not strongly 
noetherian (on either side) \cite{Ro1}.  Theorem~\ref{main theorem}(1) shows that it is also 
possible for the strong noetherian property to fail on one side only of a noetherian 
graded ring.  

It is natural to suspect that a ring for which the noetherian property fails 
after commutative base ring extension might also have strange properties when 
tensored with itself or its opposite ring.  Theorem~\ref{main theorem}(2) confirms 
such a suspicion.  The existence of a pair of finitely presented noetherian 
$k$-algebras whose tensor product is not noetherian answers \cite[Appendix, Open Problem $16'$]{GW}; our 
example shows that one can even take the algebras in question to be $\mb{N}$-graded. 

We now explain the third part of Theorem~\ref{main theorem}.  Let $A = \bigoplus_{n = 0}^{\infty} A_n$ 
be an aribtrary $\mb{N}$-graded $k$-algebra, where $k$ is an algebraically closed 
field.  In addition, assume that $A$ is \emph{connected} ($A_0 = k$) and \emph{finitely 
graded} ($\dim_k A_n < \infty$ for all $n \geq 0$).  The left noncommutative
projective scheme associated to $A$ is defined to be the pair $A \Proj 
= (A\Qgr, \mc{A})$.  Here $A\Qgr$ is the quotient category of the category of $\mb{Z}$-graded 
left $A$-modules by the full subcategory of modules which are direct limits of 
modules with finite $k$-dimension, and $\mc{A}$, called the \emph{distinguished 
object}, is the image of the module $_A A$ in $A\Qgr$.  The right noncommutative 
projective scheme $\rProj A$ of $A$ is defined analogously. The motivation for these 
definitions comes from the commutative case: if $A$ is commutative noetherian and 
$\cproj A = X$ is its associated scheme, then $A\Qgr$ and $\rQch X$ (the category of 
quasi-coherent sheaves on $X$) are equivalent categories, and $\mc{A}$ corresponds 
under this equivalence to the structure sheaf $\mc{O}_X$.

The result of Theorem~\ref{main theorem}(3) shows that noncommutative projective 
schemes associated to the two sides of a noncommutative noetherian ring may well be 
quite different.  In fact, for the ring $T$ of the theorem we will see that both  
$T\Qgr$ and $\rQgr T$ are equivalent to the category $\rQch X$ where $X = \mb{P}^d$ 
for some $d \geq 2$.  However, $\rProj T$ is isomorphic to $(\rQch X, \mc{O}_X)$, 
while $T\Proj$ is isomorphic to $(\rQch X, \mc{I})$ where $\mc{I}$ is a 
non-locally-free ideal sheaf.

Next we discuss the $\chi$ conditions, which are homological properties of graded 
rings which arose in Artin and Zhang's work in \cite{AZ94} to develop the theory 
of noncommutative projective schemes.  For each $i \geq 0$, the connected finitely graded 
$k$-algebra $A$ is said to satisfy $\chi_i$ on 
the left (right) if $\dim_k \uext^j_A(A/{A_{\geq 1}}, M) < \infty$ for all finitely 
generated left (right) $A$-modules $M$ and all $0 \leq j \leq i$, where $\uext$ 
indicates the $\ext$ group in the ungraded module category. 
If $A$ satisfies $\chi_i$ on the left for all $i \geq 0$ then we say 
that $A$ satisfies $\chi$ on the left. The $\chi_1$ condition is the most 
important of these conditions: it ensures that one can reconstruct the ring $A$ (in 
large degree) from its associated scheme $A\Proj$ \cite[Theorem 4.5]{AZ94}.  The 
other $\chi_i$ conditions for $i \geq 2$ are needed to show the finite-dimensionality 
of the cohomology groups associated to $A\Proj$ \cite[Theorem 7.4]{AZ94}.

Although the $\chi$ conditions always hold for commutative rings, Stafford and Zhang 
constructed noetherian rings for which $\chi_1$ fails on both sides \cite{StZh}.  The 
author studied rings in \cite{Ro1} which satisfy $\chi_1$ but fail $\chi_2$ on both 
sides.  Theorem~\ref{main theorem}(3) demonstrates yet more possible behaviors of the 
$\chi$ conditions:  first, that $\chi_1$ may hold on one side but not the other of a 
noetherian ring; and second, that for any $d \geq 2$ there are rings which satisfy 
$\chi_{d-1}$ but not $\chi_d$ (on one side).

Finally, we briefly describe the construction of the rings $T$ satisfying Theorem~\ref{main theorem}.
Recall that if $I$ is a left ideal in a noetherian ring $S$, then the 
\emph{idealizer} of $I$, written $\mb{I}(I)$, is the largest subring of $S$ which 
contains $I$ as a 2-sided ideal.  Explicitly, $\mb{I}(I) = \{ s 
\in S \mid Is \subseteq I \}$.  
Now let $S$ be a generic Zhang twist of a polynomial ring (see \S\ref{Zhang twist case}\ for 
the definition), which is a noncommutative graded ring generated in degree $1$.  Let 
$I$ be the left ideal of $S$ generated by a generic subspace $I_1 \subseteq S_1$ 
with $\dim I_1 = \dim S_1 -1$.  The ring $T = \mb{I}(I) \subseteq S$ is then the ring of 
interest which will satisfy properties (1)-(4) of Theorem~\ref{main theorem}.  
   
Our approach in this paper will be primarily algebraic.  Since this research was completed, the 
article \cite{KRS} has developed a geometric framework for the study of a class of algebras quite similar 
to the ones we study here.  We remark that many of the 
results below can be translated into this geometric language, which would allow one to show that 
the properties of Theorem~\ref{main theorem} hold for a wider class of idealizer rings.  Specifically, one
could work with idealizers inside twisted homogeneous coordinate rings over arbitrary integral projective schemes, 
instead of the special case of Zhang twists of polynomial rings we consider here.  Since our main purpose
is to construct some interesting examples, we will not attempt to be as general as possible and 
we will prefer the simpler algebraic constructions.

\section{Idealizer rings and the left and right noetherian property} 
\label{ideal noeth} 

As mentioned in the introduction, the main examples of this paper will be certain 
idealizer rings.  Idealizers have certainly proved useful in the 
creation of counterexamples before, but it seems that in many natural examples (for 
example those in \cite{Re} or \cite{St85}), the idealizer of a left ideal is a left 
but not right noetherian ring.  Since our intention is to create two-sided noetherian 
examples, in this brief section we will give some general characterizations of both 
the left and right noetherian properties for an idealizer ring.

Let $S$ be a noetherian ring with left ideal $I$, and let 
$T = \mb{I}(I) \subseteq S = \{s \in S \mid Is \subseteq I \}$ 
be the idealizer of $I$.  In \cite{St85}, Stafford gives a sufficient condition for 
the left noetherian property of $T$.  In the next proposition, we restate Stafford's 
result slightly to show that it characterizes the left noetherian property in case 
$S$ is a finitely generated left $T$-module, which occurs in many examples of 
interest.
\begin{proposition}
\label{char left noeth} Let $T$ be the idealizer of the left ideal $I$ of a 
noetherian ring $S$, and assume in addition that $_T S$ is finitely generated.  The 
following are equivalent:
\begin{enumerate}
\item $T$ is left noetherian.
\item $\hom_S (S/I, S/J)$ is a
noetherian left $T$-module (or $T/I$-module) for all
left ideals $J$ of $S$.
\end{enumerate}
\end{proposition}
\begin{proof}
By \cite[Lemma 1.2]{St85}, if $\hom_S (S/I, S/J)$ is a noetherian
left $T$-module for all left ideals $J$ of $S$ containing $I$,
then $T$ is left noetherian.  So if condition (2) holds, then $T$
is certainly left noetherian.

On the other hand, if $T$ is left noetherian, then since $_T S$ is
finitely generated, $_T S$ is also noetherian. Given any left
ideal $J$ of $S$, we can identify the left $T$-module $\hom_S
(S/I, S/J)$ with the subfactor $\{ x \in S \mid Ix \subseteq J
\}/J$ of $_T S$, so $\hom_S (S/I, S/J)$ is a noetherian
$T$-module. 
\end{proof}

Next, we give a characterization of the right noetherian property for idealizers of 
left ideals.  It is formally quite similar to the characterization of Proposition~\ref{char left 
noeth}, and may be of independent interest.  In 
fact, the result applies more generally to all subrings of $S$ inside of which $I$ is 
an ideal.  
\begin{proposition}
\label{right noeth} Let $S$ be a noetherian ring with left ideal $I$, and let 
$T$ be a subring of $S$ such that $I \subseteq T \subseteq \mb{I}(I)$.  The following are equivalent:
\begin{enumerate}
\item $T$ is right noetherian.
\item $T/I$ is a right noetherian 
ring, and $\tor_1^S (S/K, S/I) = (K \cap I)/KI$ is a noetherian right $T$-module
(or $T/I$-module) for all right ideals $K$ of $S$.
\end{enumerate}
\end{proposition}
\begin{proof}
The identification of $\tor_1^S (S/K, S/I)$ with the subfactor $(K \cap I)/KI$ of $T_T$
follows from \cite[Corollary 11.27(iii)]{Rotman}, and it is immediate that (1) implies (2).

Now suppose that condition (2) holds.  Since $S$ is right noetherian, $T$ is right 
noetherian if and only if $(JS \cap T)/J$ is a noetherian right $T$-module for all 
finitely generated right $T$-ideals $J$---see \cite[Lemma 6.10]{Ro1} for a proof of 
this in the graded case; the proof in the ungraded case is the same.  Let $J$ be an 
arbitrary finitely generated right ideal of $T$. Since $T/I$ is right noetherian, 
$(JS \cap T)/(JS \cap I)$ and $J/JI$ are noetherian right $T/I$-modules (the first 
injects into $T/I$, and $J$ surjects onto the second.)  Then $(JS \cap T)/J$ is right 
noetherian over $T$ if and only if 
$(JS \cap I)/JI$ is.  By \cite[Corollary 11.27(iii)]{Rotman} 
and the fact that $JSI = JI$, we may identify $(JS \cap I)/JI$ with $\tor_1^S(S/JS, 
S/I)$, which is a noetherian right module over $T$ by hypothesis.  It follows that 
$T$ is a right noetherian ring.
\end{proof}

\section{Noncommutative Proj of Graded idealizer rings} 
\label{ideal proj} Starting with this section, we focus our attention on idealizer 
rings inside connected finitely graded $k$-algebras in particular.  Our first task is 
to study the properties of the left and right noncommutative schemes associated to 
such idealizer rings, and so we begin with a review of some of the relevant 
definitions.

Below, $A$ will always be a connected finitely graded $k$-algebra, and we write 
$A\Gr$ for  the category of all $\mb{Z}$-graded left $A$-modules.  A module $M \in 
A\Gr$ is called \emph{torsion} if for every $m \in M$ there is some $n \geq 0$ such 
that $(A_{\geq n}) m = 0$.  Let $A\Tors$ be the full subcategory of $A\Gr$ consisting 
of the torsion modules, and define $A\Qgr$ to be the quotient category $A\Gr/A\Tors$, 
with quotient functor $\pi: A\Gr \to A\Qgr$.   For a $\mb{Z}$-graded $A$-module $M$ 
we define $M[n]$ for any $n \in \mb{Z}$ to be $M$ as an ungraded module, but with a 
new grading given by $M[n]_m = M_{n+m}$.  The shift functor $M \to M[1]$ is an 
autoequivalence of $A\Gr$ which naturally descends to an autoequivalence of $A\Qgr$ 
we call $s$, though we usually write $\mc{M}[n]$ instead of $s^n(\mc{M})$ for any 
$\mc{M} \in A\Qgr$ and $n \in \mb{Z}$.  

In general, any collection of data $(\mc{C}, \mc{F}, t)$ where $\mc{C}$ is an abelian 
category, $\mc{F}$ is an object of $\mc{C}$, and $t$ is an autoequivalence of 
$\mc{C}$ is called an \emph{Artin-Zhang triple}.  For every connected graded ring $A$ 
the data $(A\Qgr, \pi A, s)$ gives such a triple.  An isomorphism of two such triples 
is an equivalence of categories which commutes with the autoequivalences and under 
which the given objects correspond; see \cite[p. 237]{AZ94}.  For example, if $A$ is 
a connected graded commutative ring and $X = \cproj A$ is the associated scheme, then 
by a theorem of Serre one has that $(A\Qgr, \pi A, s)$ is isomorphic to $(\rQch X, 
\mc{O}_X, - \otimes \mc{O}(1))$.  Motivated by this, for any connected graded ring 
$A$ one calls the pair $A\Proj = (A\Qgr, \pi A)$ the \emph{left noncommutative 
projective scheme} associated to $A$, the object $\pi A$ the \emph{distinguished 
object}, and the autoequivalence $s$ of $A\Qgr$ the \emph{polarization}.  We define 
analogously the right-sided versions $\rQgr A$, $\rProj A$, etcetera of all of the 
notions above.  
  
Our analysis of the noncommutative schemes for idealizer rings will be restricted to 
rings which satisfy the following hypotheses, which will hold for a large class of 
examples we study later.  
\begin{hypothesis}
\label{base assump} Let $k$ be a field.  Let $S$ be a noetherian connected finitely $\mb{N}$-graded 
$k$-algebra, 
let $I$ be some homogeneous left 
ideal of $S$ such that $\dim_k S/I = \infty$, and put $T = \mb{I}(I)$.  Assume in 
addition that $_T S$ is a finitely generated module, and that $\dim_k T/I < \infty$.
\end{hypothesis}

Under the assumptions of \ref{base assump}, both the left and right 
noncommutative schemes for the 
idealizer ring $T$ are closely related to those for the ring $S$, as we see now.  
\begin{lemma}
\label{cat equivs} Assume Hypothesis~\ref{base assump}.
\begin{enumerate}
\item There is an isomorphism of triples
$(S\Qgr, \pi I , s) \cong (T\Qgr, \pi T , s)$.
\item There is an isomorphism of triples
$(\rQgr S, \pi S, s) \cong (\rQgr T, \pi T , s)$.
\end{enumerate}
\end{lemma}
\begin{proof}

(1) 
Suppose that $M \in S\Gr$.  Then we claim 
that if $_T M \in T\Tors$, then $_S M \in S\Tors$.  To prove this fact, note first that 
if $_T M$ is finitely generated, then $M$
if finite-dimensional over $k$, so obviously $_S M \in S\Tors$.  In general, $_T M$
is a direct limit of finite-dimensional $T$-modules, so $M' = S \otimes_T M$ is a direct limit 
of finite-dimensional $S$-modules and 
thus $M' \in S\Tors$.  Since there is an $S$-module surjection $M' \to M$, this
completes the proof of the claim.  

Now we define two functors by the rules
\begin{eqnarray*}
F: T\Gr & \lra & S\Gr \\
    _T M & \mapsto & _S(I \otimes_T M) \\
G: S\Gr & \lra & T\Gr \\
    _S N & \mapsto & _T N
\end{eqnarray*}
together with the obvious actions on morphisms.
If $_T M \in T\Gr$, then since $\dim_k T/I < \infty$ it follows by calculating using 
a free resolution of $M$ that $\utor_j^T(T/I, M)$ is a torsion left $T$-module for 
all $j \geq 0$.  Then the natural map $I \otimes_T M \ra T \otimes_T M = M$ has 
torsion kernel and cokernel for all $M \in T\Gr$.  In particular, if $M \in T\Tors$ 
then $F(M) \in T\Tors$, so $F(M) \in S\Tors$ by the earlier claim.  It 
follows that $F'= \pi \circ F: T\Gr \to S\Qgr$ is an exact functor, and that 
$F'(M) = 0$ for all $M \in T\Tors$.  Then by the universal property 
of the quotient category \cite[Corollary 4.3.11]{Po}, 
$F'$ descends to a functor $\overline{F}: T\Qgr \to S\Qgr$.
Similarly, it is clear that if $N \in S\Tors$ then $G(N) = N \in T\Tors$.  Then 
$G' = \pi \circ G: S\Gr \to T\Qgr$
is an exact functor with $G'(N) = 0$ for all $N \in S\Tors$, so $G'$
descends to a functor $\overline{G}: S\Qgr \to T\Qgr$.

For $M \in T\Gr$, since the natural multiplication 
map $GF(M) = I \otimes M \ra M$ has torsion kernel and cokernel, it descends to an isomorphism 
$\overline{G}\overline{F}(\pi M ) \cong \pi M$.  If $N \in S\Gr$, then the 
multiplication map $FG(N) = I \otimes_T N \ra N$ is actually an $S$-module map.  The kernel and 
cokernel of this map are torsion $T$-modules, and therefore torsion $S$-modules by 
the claim above.  Thus we also get a isomorphism $\overline{F}\overline{G}(\pi N) \cong \pi N$ in $S\Qgr$.

We conclude that $\overline{F}$ and $\overline{G}$ are inverse equivalences of 
categories.  Moreover, obviously $\overline{F}(\pi T) \cong \pi I$, and all of the 
maps are compatible with the shift functors $s$, since $F$ and $G$ are compatible 
with the shift functors in the categories $S\Gr$ and $T\Gr$.

(2) Because $SI = I \subseteq T$, we have $(S/T)I = 0$ and so since $T/I$ is finite
dimensional we see that $(S/T)_T$ is torsion.  By assumption we also know that $_T 
(S/T)$ is finitely generated.  Now the proof of this triple isomorphism is entirely 
analogous to the proof of \cite[Proposition 2.7]{StZh}, with the exception that the 
authors assume there that $T$ is noetherian and then prove the required equivalence 
for the subcategories of noetherian objects.  
We leave it to the reader to make the obvious adjustments to the proof to show 
without the noetherian assumption 
that $(\rQgr S, \pi S, s) \cong (\rQgr T, \pi T, s)$.
\end{proof}

\begin{remark}
The graded idealizer rings studied by Stafford and Zhang in \cite{StZh} have the 
special property that the ideal $I$ is a principal ideal generated by an element of degree $1$ 
in a graded Goldie domain $S$.  In 
that case, $T = \mb{I}(I)$ is isomorphic to its opposite ring, and 
thus the differences between parts (1) and (2) of Proposition~\ref{cat equivs} must 
disappear (indeed, in this case $\pi I \cong \pi S[-1]$).   In the general case, 
however, it is clear from Proposition~\ref{cat equivs} that we should expect the 
noncommutative schemes $T\Proj$ and $\rProj T$ to be non-isomorphic.  
\end{remark}

The information provided by Lemma~\ref{cat equivs} will allow us to prove
with ease several further results about the  
noncommutative projective schemes of idealizer rings.   
First, we may show in wide generality that passing to a Veronese ring of $T$ does not 
affect the associated noncommutative projective schemes.  Recall that for an 
$\mb{N}$-graded ring $A$ the $n$th Veronese ring of $A$ is the graded ring 
$A^{(n)} = \bigoplus_{i = 0}^{\infty} A_{in}$.
\begin{proposition}
\label{Veroneses} Assume Hypothesis~\ref{base assump}, and in addition let $S$ be 
generated in degree $1$.  Choose $n \geq 1$ and write $T' = T^{(n)}$, $S' = S^{(n)}$, 
and $I' = I^{(n)} = \bigoplus_{i = 0}^{\infty} I_{in}$.  Let $R' \subseteq S'$ be the idealizer of 
the left ideal $I'$ of $S'$.       
\begin{enumerate}
\item $T'$ and $R'$ are isomorphic in large degree. 
\item There are isomorphisms of 
noncommutative projective schemes $T\Proj \cong T'\Proj$ and $\rProj T \cong \rProj T'$.
\end{enumerate}
\end{proposition}
\begin{proof}
(1) As ungraded rings, we may identify $R'$, $T'$ and $S'$ with subrings of $S$.
Suppose that $x \in (R'_m)$, so that 
$I' x \subseteq I'$.  Then since the left ideal $I$ of $S$ is generated in 
some finite degree, we see that in 
the ring $S$ we have $(I_{\geq p}) x \subseteq I$ for some $p \geq 0$, where 
$x \in S_{nm}$.  Since the torsion submodule of $_S(S/I)$ is finite dimensional, if $m 
\gg 0$ then $I x \subseteq I$ and hence $x \in T$.  Then as an element of $S'$, $x \in T'$.  Since
the inclusion $T' \subseteq R'$ is obvious, $T'$ and $R'$ must agree in large degree.

(2) Since $_T S$ is finitely generated and $\dim_k T/I < \infty$, we see that $_{T'} S'$ is 
finitely generated and $\dim_k T'/I' < \infty$.  Then because 
$T'$ and $R'$ agree in large degree by part (1), it follows that $_{R'} 
S'$ is finitely generated and that $\dim_k R'/I' < \infty$.  Also, since $S$ is 
noetherian, $S'$ must be noetherian \cite[Proposition 5.10(1)]{AZ94}.  

Now we claim that we have isomorphisms of noncommutative projective schemes 
\[
(\rQgr T, \pi T) \cong (\rQgr S, \pi I) \cong (\rQgr S', \pi I') \cong 
\rQgr (R', \pi R') \cong \rQgr (T', \pi T').
\]
To see this, note that since $S$ is generated in degree $1$, there is an isomorphism 
$\rProj S \cong \rProj S'$ \cite[Proposition 5.10(3)]{AZ94}; the associated 
equivalence of categories  $\rQgr S \simeq \rQgr S'$ sends $\pi I$ to $\pi 
I'$.  The second isomorphism follows, and the first and third follow from 
Proposition~\ref{cat equivs}(1), applied to $T \subseteq S$ and to $R'\subseteq 
S'$, respectively.  Last, the final isomorphism follows from part (1).  
Altogether this chain of isomorphisms says that $\rProj T \cong \rProj T'$.  

The argument on the left side is very similar, except using the other triple 
isomorphism of Proposition~\ref{cat equivs}, and is left to the reader.  
\end{proof}

Next we will show that 
under mild hypotheses the noncommutative projective schemes 
associated to $S$ and $T$ (on either side) have the same cohomological dimension; we review
the definition of this property now.   
Cohomology groups for the noncommutative projective scheme $A\Proj$ are defined 
by setting $\coH^i(\mc{M}) = \ext^i_{A\Qgr}(\pi A, \mc{M})$ for all $\mc{M} \in A\Qgr$.  Then 
the \emph{cohomological dimension} of $A\Proj$ is  
\[
\cd(A\Proj) = \max \{i \mid \coH^i(\mc{M}) \neq 0\ \text{for some}\ \mc{M} \in A\Qgr 
\}
\]
and the \emph{global dimension} of the category $A\Qgr$ is 
\[
\gd(A\Qgr) = \max \{i \mid \ext^i_{A\Qgr}(\mc{M}, \mc{N}) \neq 0\ \text{for some}\ 
\mc{M}, \mc{N} \in A\Qgr \}.
\]
The right-sided versions of these notions are defined similarly.

\begin{proposition}
\label{cd} 
Assume Hypothesis~\ref{base assump}.
\begin{enumerate}
\item $\cd(\rProj T) = \cd(\rProj S)$.
\item Assume in addition that $S$ is a domain with $\gd(S\Qgr) = \cd(S\Proj) < \infty$.
Then $\cd(T\Proj) = \cd(S\Proj)$.
\end{enumerate}
\end{proposition}
\begin{proof}
(1) This part is immediate from the triple isomorphism of Proposition~\ref{cat 
equivs}(2). 

(2) By proposition~\ref{cat equivs}(1), we have the isomorphism of triples $(T\Qgr, \pi 
T, s) \cong (S\Qgr, \pi I, s)$.  From this it quickly follows that 
\[
\cd(T\Proj) \leq \gd(T\Qgr) = \gd(S\Qgr) = \cd(S\Proj). 
\]  
Let $d = \cd(S\Proj)$.  To finish the proof that $\cd(T\Proj) = \cd(S\Proj)$ we have only to show that 
there is some $\mc{F} \in S\Qgr$ such that $\ext^d_{S\Qgr}(\pi I, \mc{F}) \neq 0$.  
Since $S$ is a domain, we may choose some injection $S[-m] \to I$ for some $m \geq 
0$, and passing to $S\Qgr$ we have a short exact sequence $0 \to \pi S[-m] \to \pi I 
\to \mc{N} \to 0$ for some $\mc{N}$.  Since $S\Proj$ has cohomological dimension $d$ 
we may choose some $\mc{F} \in S\Qgr$ with $\ext^d_{S\Qgr}(\pi S[-m], \mc{F}) \neq 
0$.  But $\ext^{d+1}_{S\Qgr}(\mc{N}, \mc{F}) = 0$ since the global dimension of 
$S\Qgr$ is $d$, so we conclude from the long exact sequence in $\ext$ that $\ext^d_{S\Qgr}(\pi 
I, \mc{F}) \neq 0$.
\end{proof}

\section{The $\chi$ conditions for graded idealizers}
\label{ideal chi} 

The goal of this section is to begin an analysis of the $\chi$ conditions, which we 
defined in the introduction, for the case of graded idealizer rings $T$ satisfying 
Hypothesis~\ref{base assump}.  The main result below  will show that if $S$ itself 
satisfies left $\chi$, then the left $\chi$ conditions for the idealizer ring $T$ may 
be characterized in terms of homological algebra over $S$ only.  
We also study the right $\chi$ conditions for $T$;  
the analysis of these turns out to be a much simpler matter.  

We review several definitions which we will need before 
proving the main result of this section.  
A module $M \in A\Gr$ is \emph{right bounded} if $M_n = 0$ for $n \gg 0$, \emph{left 
bounded} if $M_n = 0$ for $n \ll 0$, and \emph{bounded} if it is both left and right 
bounded.  $M$ is \emph{finitely graded} if $\dim_k M_n <\infty$ for all $n \in 
\mb{Z}$.  For $M,N \in A\Gr$, $\hom_A(M,N)$ means the group of degree-preserving 
module homomorphisms, and $\ext_A^i(M, -)$ is the $i$th right derived functor of $\hom_A(M, -)$.  
We also set $\uhom_A(M,N) = \bigoplus_{n \in \mb{Z}} \hom(M, 
N[n])$, which is the same as the group of homomorphisms in the ungraded category if 
$M$ is finitely generated. More generally, we write 
$\uext_A^i(M,N) = \bigoplus_{n \in \mb{Z}} \ext_A^i(M, N[n])$.  We make 
similar definitions in the category $A\Qgr$; so 
$\uext_{A\Qgr}^i(\mc{M}, \mc{N}) = \bigoplus_{n \in \mb{Z}} \ext_{A\Qgr}^i(\mc{M}, \mc{N}[n])$.  
Finally, let $A\gr$ be the subcategory of all noetherian modules in $A\Gr$.

Note that we have defined the $\chi$ conditions for not necessarily noetherian 
algebras; it is easy to prove, however, that the left $\chi_0$ condition for a 
connected graded ring $A$ is equivalent to the left noetherian property for $A$.  
Recall also that if $A$ is connected graded left noetherian with modules $M \in A\gr$ 
and $N \in A\Gr$, then for any $j \geq 0$ we have $\uext^j_{A\Qgr}(\pi M, \pi N) 
\cong \lm \uext^j_A(M_{\geq n}, N)$ \cite[Proposition 7.2(1)]{AZ94}.  In particular, 
in this case there is a natural map of vector spaces $\uext^j_A(M, N) \to 
\uext^j_{A\Qgr}(\pi M, \pi N)$.  In the proof of the following proposition we will 
use several results of Artin and Zhang from \cite{AZ94} which interpret the $\chi$ 
conditions in terms of the properties of such maps.

\begin{proposition}
\label{char chi} Assume Hypothesis~\ref{base assump}, and assume also 
that $S$ satisfies $\chi$ on the left.  Then $T$ 
satisfies $\chi_i$ on the left for some $i \geq 0$ if and only 
if $\dim_k \uext^j_S(S/I,M) < \infty$ for all $0 \leq j \leq i$ and all $M \in S\gr$.
\end{proposition}
\begin{proof}
Since any $M \in S\gr$ has a finite filtration by cyclic $S$-modules, it follows from 
Theorem~\ref{char left noeth} that $T$ is left noetherian if and only if 
$\uhom_S(S/I,M)$ is a noetherian left $T/I$-module (equivalently, of finite 
$k$-dimension) for all $M \in S\gr$.  Since the left noetherian property for $T$ is 
equivalent to left $\chi_0$ for $T$ (as we remarked before the proposition), the 
characterization of the proposition holds when $i = 0$.

Now assume that $T$ is left noetherian.  There is an isomorphism of triples $(S\Qgr, 
\pi I, s) \cong (T\Qgr, \pi T, s)$ by Proposition~\ref{cat equivs}(1).  For any $M 
\in S\gr$ we have a diagram
\[
\xymatrix@C=2.2em{ M \ar[r]^{\gamma} \ar[d]^{\alpha} & \uhom_S(I, M) \ar[d]^{\beta} \\
\uhom_{T\Qgr}(\pi T, \pi M) \ar[r]^{\cong} & \uhom_{S\Qgr}(\pi I, \pi M) }
\]
where the bottom arrow is an isomorphism by the triple isomorphism, $\alpha$ and 
$\beta$ are the natural maps, and $\gamma$ is part of the long exact sequence in 
$\uext$.  It is straightforward to check that this diagram commutes.  Now since $S$ 
has $\chi$, the map $\beta$ is an isomorphism in large degree \cite[Proposition 
3.5(3)]{AZ94}.  Furthermore, $\chi_1$ holds on the left for $T$ if and only if the 
map $\alpha$ has right bounded cokernel for all $M \in T\gr$ \cite[Proposition 
3.14(2a)]{AZ94}. Note that it is equivalent to require that $\alpha$ have bounded cokernel
for all $M \in S\gr$, as follows:  if $M \in T\gr$, then $IM \in S\gr$ with 
$\dim_k M/IM < \infty$ and thus $\pi M = \pi IM$; conversely, if $M \in S\gr$ 
then $M \in T\gr$ since $_T S$ is 
finitely generated and $T$ is left noetherian.  Thus from the diagram it follows that 
$\chi_1$ holds for $T$ on the left if and only if $\gamma$ has right bounded cokernel 
for all $M \in S\gr$.  But the cokernel of $\gamma$ is $\uext^1_S(S/I, M)$, which is 
always finitely graded and left bounded, so is right bounded if and only if it has 
finite $k$-dimension.  Thus the proposition holds for $i = 1$. 

Next, assume that $\chi_1$ holds on the left for $T$.  Then the proof of the 
noncommutative version of Serre's finiteness theorem \cite[Theorem 7.4]{AZ94} 
shows that $\chi_i$ holds for $T$ for some $i \geq 2$ if and only if for every $M \in 
T\gr$, the graded cohomology group $\ucoH^j(\pi M) =  \uext^j_{T\Qgr}(\pi T, \pi M)$ 
is finitely graded for all $0 \leq j < i$ and right bounded for all $1 \leq j < i$.  
Similarly as in the last paragraph, one sees that it is equivalent to require this 
condition for all $M \in S\gr$.  Now for every $M \in S\gr$ and $j \geq 1$ we have a 
sequence of maps 
\[
\xymatrix@C=1.5em{ \uext^{j+1}_S(S/I, M) \ar[r]^-{\cong} & \uext^j_S(I, M) \ar[r]^-{\alpha} & 
\uext^j_{S\Qgr}(\pi I, \pi M) \ar[r]^-{\cong} & \uext^j_{T\Qgr}(\pi T, \pi M) }
\]
where the first isomorphism comes from the long exact sequence in $\uext$, the 
natural map $\alpha$ is an isomorphism in large degree since $S$ satisfies $\chi$ 
\cite[Proposition 3.5(3)]{AZ94}, and the final isomorphism comes from the isomorphism 
of triples in Proposition~\ref{cat equivs}(1).  In addition, 
$\uext^j_{S\Qgr}(\pi I, \pi M)$ is always finitely graded 
for any $j$, since $S$ has $\chi$ \cite[Corollary 7.3(3)]{AZ94}.  Thus we see 
altogether that, assuming $\chi_1$ holds for $T$, $\chi_i$ holds for $T$ for some $i 
\geq 2$ if and only if $\uext^j_S(S/I, M)$ is right bounded (equivalently, finite 
dimensional over $k$ since it is always left bounded and finitely graded) for all $2 
\leq j \leq i$ and all $M \in S\gr$.  This 
proves the characterization of $\chi_i$ for $i \geq 2$, and 
concludes the proof of the proposition.
\end{proof}

In contrast to Proposition~\ref{char chi}, on the right side only the $\chi_0$ 
condition for $T$ (equivalently, the right noetherian property for $T$) is 
potentially subtle to analyze.  The higher $\chi$ conditions automatically must fail, 
as follows.
\begin{proposition}
\label{right chi}
Assume Hypothesis~\ref{base assump}.  Then $T$ fails $\chi_i$ on the right for all $i \geq 1$.
\end{proposition}
\begin{proof}
We may assume that $T$ is right noetherian, since otherwise right $\chi_0$ fails for 
$T$ and so by definition right $\chi_i$ fails for all $i \geq 0$.  Also, we 
need only show that $T$ fails right $\chi_1$.  For this, 
the same argument outlined in \cite[p. 424]{StZh} works here; since it is simple we 
briefly repeat it.  By hypothesis, we have $SI = I$, $\dim_k T/I < \infty$, and 
$\dim_k S/I = \infty$.  So the natural map 
\[
T \to \uhom_{\rQgr T}(\pi T, \pi T) = \uhom_{\rQgr T}(\pi I, \pi I)
\]
has a cokernel which is not right bounded, since $S \subseteq \uhom_{\rQgr T}(\pi I, \pi I)$.  
Then by \cite[Proposition 3.14(2a)]{AZ94}, $T$ must fail $\chi_1$ on the right.
\end{proof}

\section{Idealizers inside Zhang twists of polynomial rings}
\label{Zhang twist case} 
In the current section, we introduce a special class of 
graded idealizers on which we will focus for the remainder of the paper. 

Fix a commutative polynomial ring $U = k[x_0, x_1, \dots, x_d]$ in $d+1$ variables, 
and some graded automorphism $\phi$ of $U$. 
Let $S$ be the \emph{left Zhang twist} of $U$ by $\phi$.  This is a new ring which 
has the same underlying $k$-space as the ring $U$, but a new multiplication defined 
by the rule $fg = \phi^n(f) \circ g$ for $f \in S_m, g \in S_n$, where $\circ$ is the 
multiplication in $U$.  We continue this same notational convention throughout, 
whereby juxtaposition means multiplication in $S$ and the symbol $\circ$ appears when 
the commutative multiplication in $U$ is intended.

It is immediate that $S$ is a noetherian domain \cite[Theorem 1.3]{Zhang}. One may 
also twist modules:  given a graded $U$-module $M$, one may form a graded left 
$S$-module with the same underlying vector space as $M$ but with $S$-action $fg = 
\phi^n(f) \circ g$ for $f \in S_m, g \in M_n$, where again $\circ$ indicates the 
$U$-action.  In this way we get a functor $U\Gr \to S\Gr$ which is an equivalence of 
categories \cite[Corollary 4.4(1)]{Zhang}.  In particular, the graded left ideals of 
$S$ and the graded (left) ideals of $U$ are in one-to-one correspondence, and if $J$ 
is a graded left $S$-ideal we use the same name $J$ for the corresponding graded 
$U$-ideal.

Now we will idealize left ideals of $S$ which are generated by a codimension-$1$ 
subspace of the elements of degree $1$. 
Specifically, from now on we will consider the following hypothesis and notations.  
\begin{hypothesis}
\label{second assump} Let $k$ be an algebraically closed base field.  Choose some $d 
\geq 2$, a point $c \in \mb{P}^d$, and an automorphism $\varphi \in \aut \mb{P}^d$.  
Let $\phi$ be a graded automorphism of $U = k[x_0, \dots, x_d]$ such that $\varphi$ 
is the corresponding automorphism of $\cproj U = \mb{P}^d$, and define $S = S(\varphi)$ 
to be the left Zhang twist of $U = k[x_0, \dots, x_d]$ by the automorphism $\phi$.  
(Although the automorphism $\phi$ corresponding to $\varphi$ is determined only up to scalar 
multiple \cite[Example 7.1.1]{H}, it is easy to check that changing $\phi$ 
by a nonzero scalar does not change the ring $S$ up to isomorphism.)  Let 
$I$ be the left ideal of $S$ consisting of all homogeneous elements vanishing at the 
point $c$.  Define $T = T(\varphi, c) = \mb{I}(I) \subseteq S$. Also write $c_n = 
\varphi^{-n}(c)$ for $n \in \mb{Z}$. 
\end{hypothesis}

In general, the properties of the 
ring $T = T(\varphi, c)$ depend on the properties of the orbit $\mc{C} = \{ c_n \}_{n \in \mb{Z}}$.  
We are most interested in the ``generic'' case, and 
so we will usually assume at least that $\mc{C}$ is 
infinite.  Under such an assumption, we see next 
that the idealizer rings $T$ have the following basic properties.
\begin{lemma}
\label{basic props}
Assume Hypothesis~\ref{second assump}.  If the points $\{ c_n \}_{n \in \mb{Z}}$ are all 
distinct, then 
\begin{enumerate}
\item $T = k + I$.
\item $T^{(n)}$ is not generated in degree $1$ for any $n \geq 1$. 
\item $\dim_k (S/IS) < \infty$.
\item $_T S$ is finitely generated.
\item $T$ is a finitely generated $k$-algebra.
\item Hypothesis~\ref{base assump} is satisfied.
\end{enumerate}
\end{lemma}
\begin{proof}
(1) 
We have $T_n = \{x \in S_n \mid Ix \subseteq I\}$.  If $\phi^n(I) \neq I$, then since 
$I$ is prime in $U$, $\phi^n(I) \circ x \subseteq I$ forces $x \in I$.  Since we assume
that $c$ has infinite order under $\varphi$, $\phi^n(I) \neq I$ for all $n \neq 0$ 
and so $T_n = I_n$ for $n \geq 1$.

(2) 
If $T^{(n)}$ were generated in degree one for some $n \geq 1$, then would we have $T_n T_n = 
T_{2n}$, which in the commutative ring $U$ translates to $\phi^n(I)_n \circ I_n = 
I_{2n}$.   Since $I$ and $\phi^n(I)$ are different homogeneous prime ideals of $U$ 
which are generated in degree $1$, it is easy to see that such an equation is 
impossible.  

(3) Set $J = IS$.  We have that $J = \sum_{i = 0}^{\infty} I S_i = \sum_{i = 0}^{\infty} 
\phi^i(I) \circ U_i$.  Since the points $\{c_i\}$ are all distinct, it is clear 
that the vanishing set of the ideal $J$ in $\mb{P}^d$ is empty.  Thus $\dim_k U/J < \infty$ by the 
graded Nullstellensatz; equivalently, $\dim_k (S/IS) < \infty$.

(4) By the graded Nakayama lemma, a $k$-basis of $S/T_{\geq 1}S = S/IS$ is a minimal
generating set for $_T S$, so (4) follows immediately from (3).

(5) $T$ is generated as a $k$-algebra by some elements $t_i \in T_{\geq 1}$ if and 
only if $T_{\geq 1}$ is generated as left $T$-ideal by the $t_i$; so to prove (5) we just need 
to show that $_T I$ is finitely generated.  Since by 
part (4) we know that $_T S$ is finitely generated, we have $T S_{\leq n} = S$ for some $n \geq 0$.  Then 
$T S_{\leq n} T_1 = S T_1 = I$ is a finitely generated left 
$T$-module. 

(6) Since $\dim_k I_n = \dim_k S_n - 1$ for all $n \geq 1$, it is clear that $\dim_k S/I = \infty$.  
The other necessary properties follow from (1) and (4).
\end{proof} 

The noetherian property on the left is also straightforward to analyze.
\begin{proposition}  
\label{left noeth} Assume Hypothesis~\ref{second assump}, and that the points $\{ c_n \}_{n \in \mb{Z}}$
are all distinct.  Then $T$ is left noetherian.
\end{proposition}
\begin{proof}
We have that $T = k + I$ and that $_T S$ is finitely generated, by 
Lemma~\ref{basic props}.  Thus the hypotheses of 
Theorem~\ref{char left noeth} are satisfied and to show that 
$T$ is left noetherian we need to show that 
$\uhom_S(S/I, S/J)$ is a left noetherian (equivalently, finite dimensional) $T/I = 
k$-module for all graded left ideals $J$ of $S$.  Using the equivalence of categories 
$S\Gr \sim U\Gr$ and the existence of prime filtrations in $U$, we see that every 
cyclic graded left $S$-module $S/J$ has a finite graded filtration with factors of the form $S/L$ 
where $L$ is prime when considered as an ideal of $U$.  Thus we may reduce to the case
that $J$ is a prime ideal of $U$.  If $J = U_{\geq 1}$ then 
obviously $\uhom_S(S/I, S/J)$ is finite dimensional, so we also 
may assume that $J \neq U_{\geq 1}$.

Now we may make the identification of vector spaces 
\[ \uhom_S(S/I, S/J)_n = \{ x \in U_n \mid \phi^n(I) \circ x \subseteq J \}/J_n.  
\]
Since 
the points $\{c_i \}_{i \geq 0}$ are distinct, $\phi^n(I) \subseteq J$ can occur for at 
most one value of $n$; since $J$ is prime, we see that 
$\{ x \in U_n \mid \phi^n(I) \circ x \subseteq J \} = J_n$ for all $n \gg 0$ and so 
$\uhom_S(S/I, S/J)_n = 0$ for $n \gg 0$.  Thus 
$\uhom_S(S/I, S/J)$ is indeed finite dimensional over $k$.
\end{proof}

The right noetherian property and the left $\chi$ conditions for the ring $T$ 
depend on a more subtle property of the set of points $\{c_n \}_{n \in \mb{Z}}$.  Given a 
subset $\mc{C}$ of closed points 
of $\mb{P}^d$, we say that $\mc{C}$ is \emph{critically dense} if 
every infinite subset of $\mc{C}$ has Zariski closure equal to all of $\mb{P}^d$.
\begin{proposition}
\label{main noeth and chi result} Assume Hypothesis~\ref{second assump}, and assume
in addition that the set of points $\{c_n\}_{n \in \mb{Z}}$ is critically dense in $\mb{P}^d$.  Then 
\begin{enumerate}
\item $T$ satisfies left $\chi_{d-1}$ but fails left $\chi_d$. 
\item $T$ is right noetherian.
\end{enumerate}
\end{proposition}
\begin{proof}
(1) By \cite[Lemma 8.4(2)]{Ro1}, if $J$ is a graded left ideal of $S$ then we have 
\[
\uext^i_S(S/I,S/J)_n \cong \uext^i_U (U/I, U/\phi^{-n}(J))_n
\]
as $k$-spaces, for each $n \in \mb{Z}$.  It follows that $\uext^d_S(S/I, S)_n \cong \uext^d_U(U/I, U)_n \neq 0$
for all $n \geq 0$, since one may calculate that $\uext^d(U/I, U) \cong (U/I)[d]$ easily
from a Koszul resolution of $U/I$.  So $S$ fails $\chi_d$ on the left.  On the other hand, 
\cite[Proposition 8.6(1)]{Ro1} proves that since $\{c_n \}_{n \in \mb{Z}}$ is critically dense,
we have $\dim_k \uext^i_S(S/I, M) < \infty$ for all $0 \leq i \leq d-1$ and all finitely generated
left $S$-modules $M$.  Then $T$ satisfies $\chi_{d-1}$ on the left by Proposition~\ref{char chi}.

(2)
If we can show that every module of the form
$(JS \cap T)/J$, for $J$ a finitely generated right ideal of $T$, is finite dimensional, then 
\cite[Lemma 5.10]{Ro1} shows that $T$ is right noetherian.  Note that $T$ is an Ore
domain, since it is a domain of finite GK-dimension \cite[Proposition 4.13]{KL}.  
Then the same proof as in \cite[Lemma 5.9]{Ro1} shows that every module of the form
$(JS \cap T)/J$ for $J$ a finitely generated right ideal of $T$ is 
filtered by subfactors of modules of the form $(fS \cap T)/fT$ and $S/(fS + T)$ for nonzero 
homogeneous $f \in T$.
Thus we will just need to prove that modules of those forms are finite dimensional over $k$.

Recall that $T = k + I$ by Lemma~\ref{basic props}(1).  
Fix $n \geq 1$ and let $f \in T_n$ be arbitrary.  We have for $m \geq n$ that 
$(fS + T)_m = \phi^{m-n}(f) \circ U_{m -n} + I_m$.  Since $T_m = I_m$ has codimension 
$1$ inside $S_m$ for all $m \geq 1$ and $I$ is prime in $U$, this implies that $(fS + T)_m = 
S_m$ if and only if $\phi^{m-n}(f) \not \in I$.  

Similarly, again assuming $m \geq n$, we have $(fS \cap T)_m = (\phi^{m-n}(f) \circ 
U_{m -n}) \cap I_m$.  If $\phi^{m-n}(f) \not \in I$, then as $I$ is prime, 
$(\phi^{m-n}(f) \circ U_{m -n}) \cap I_m = \phi^{m-n}(f) \circ I_{m -n} = (fT)_m$.  Conversely, if 
$\phi^{m-n}(f) \in I$, then $(fS \cap T)_m = (fS)_m \neq (fT)_m$. 

Now since $\{ c_n \}_{n \in \mb{Z}}$ is a critically dense set of points, every 
homogeneous $f \in S$ satisfies $f \not \in \phi^n(I)$ for $n \ll 0$, which is 
equivalent to $\phi^n(f) \not \in I$ for $n \gg 0$.  We conclude that for any homogeneous 
$0 \neq f \in T$ the modules $(fS \cap T)/fT$ and $S/(fS + T)$ are finite dimensional, as required.
\end{proof}

\section{The strong noetherian property}
\label{tensors}  We continue to study idealizer rings $T$ satisfying 
Hypothesis~\ref{second assump},
and we maintain the notation introduced in the previous section.  In \cite{Ro1}, the 
author showed the existence of rings which are not strongly noetherian on either 
side.  Here we will show that the idealizer rings $T$ are typically strongly 
noetherian on one side but not the other.

Let $A$ be an arbitrary $k$-algebra.
We call a left $A$-module $M$ \emph{strongly noetherian} if $M \otimes_k B$ is 
a noetherian left $A \otimes_k B$-module for every commutative noetherian $k$-algebra $B$.  More generally,
$M$ is \emph{universally noetherian} if $M \otimes_k B$ is noetherian over $A \otimes_k B$ for 
every noetherian $k$-algebra $B$.
\begin{proposition}
\label{strong noeth} Assume Hypothesis~\ref{second assump}, and assume further that the set of points 
$\{ c_n \}_{n \in \mb{Z}}$ is critically dense.  Then $T$ is a noetherian ring such that 
\begin{enumerate}
\item $T$ is universally left noetherian.
\item $T$ is not strongly right noetherian.
\end{enumerate}
\end{proposition}
\begin{proof} 
That $T$ is noetherian follows from Propositions~\ref{left noeth} and 
\ref{main noeth and chi result}.  

(1) 
We note that the ring $S$ is universally left noetherian, as follows.  For 
any noetherian $k$-algebra $B$, the ring $U \otimes_k B \cong B[x_0, \dots, x_d]$ is 
noetherian by the Hilbert basis theorem.  Then since $S \otimes_k B$ is a left Zhang twist 
of $U \otimes_k B$, it is also left noetherian \cite[Theorem 1.3]{Zhang}. 
Now we prove that $T$ is universally noetherian on the left.  We know that $M = {}_T 
(S/T)$ is finitely generated by Lemma~\ref{basic props}, and since $\dim_k M_n = 1$ 
for all $n \geq 1$ we see that $M$ must have Krull dimension $1$.  By 
\cite[Theorem 4.23]{ASZ}, $M$ is a universally noetherian left $T$-module.  So 
if $B$ is any noetherian $k$-algebra, then $M \otimes_k B = 
(S \otimes_k B)/(T \otimes_k B)$ is a noetherian left $T \otimes_k B$-module.  Then by 
\cite[Lemma 4.2]{ASZ}, since $S \otimes_k B$ is left 
noetherian, $T \otimes_k B$ is also left noetherian.

(2) 
The proof which we now present that $T$ is not strongly noetherian on the right is quite 
analogous to the proof in \cite[\S 7]{Ro1} that the ring $R$ studied in that paper 
is not strongly noetherian.  Let us first make a few comments about notation.  We  
use subscripts to indicate extension of scalars, for example $U_B = U \otimes_k B$.  The automorphism 
$\phi$ of $U$ naturally extends to an automorphism of $U_B$ such that $S_B$ is 
again the left Zhang twist of $U_B$ by $\phi$.  
We extend also our notational convention, so that juxtaposition 
means multiplication in $S_B$ and $\circ$ means the commutative multiplication in 
$U_B$.  Fix once and for all some particular choice of homogeneous coordinates for each of 
the points in $\{c_n\}_{n \in \mb{Z}} \subseteq \mb{P}^d_k$.  Then for $f \in U_B$, 
the expression $f(c_n)$ denotes polynomial evaluation at the fixed coordinates for 
$c_n$, giving a well-defined value in the ring $B$.  

Because by assumption the point set $\{c_n\}_{n \in \mb{Z}}$ is critically dense, the 
same proof as in \cite[Theorem 7.4]{Ro1} shows that there exists a noetherian commutative 
$k$-algebra $B$ which is a unique factorization domain, constructed as an 
\emph{infinite affine blowup} of affine space, and containing elements $f, g \in (U_B)_1$ with 
the following properties:
\begin{enumerate}
\item $g(c_i) = \Omega_i f(c_i)$ for some $\Omega_i \in B$, for all $i
\leq 0$. 
\item For all $i \ll 0$, $f(c_i)$ is not a unit in $B$.
\item $\gcd(f, g) = 1$ in $U_B$.
\end{enumerate}

Note that a homogeneous element $f \in U_B$ is in $(T_B)_{\geq 1} = I \otimes_k B$ if and only if $f(c_0) = 0$.  Now 
for each $n \geq 1$ we may choose some element  
$\theta_n \in (S_B)_n \setminus (T_B)_n$ with coefficients in $k$. Putting $t_n = 
(\Omega_{-n} f -g) \theta_n$, we have in terms of the commutative multiplication in 
$U_B$ that $t_n = \phi^n(\Omega_{-n} f - g) \circ \theta_n$, and since 
$\phi^n(\Omega_{-n} f -g)(c_0) = (\Omega_{-n} f -g)(c_{-n}) = 0$ we see that $t_n \in 
(T_B)_{n+1}$.  Suppose for some $n$ that $t_{n+1} = \sum_{i = 1}^n t_i r_i$ with $r_i 
\in (T_B)_{n-i +1}$.  Then
\[
\phi^{n+1}(\Omega_{-n-1} f -g) \circ \theta_{n+1} = \sum_{i = 1}^n
\phi^{n+1}(\Omega_{-i} f - g) \circ \phi^{n-i +1}(\theta_i) \circ r_i.
\]

Rewriting this equation in the form $h_1 \circ \phi^{n+1}(f) = h_2
\circ \phi^{n+1}(g)$, and using that $\gcd(f, g) =1$, we may
conclude that $\phi^{n+1}(g)$ divides $h_1$, where
\[
h_1 = \Omega_{-n-1} \theta_{n+1} - \sum_{i =1}^n \Omega_{-i} \phi^{n-i 
+1}(\theta_i) \circ r_i.
\]
Then $(\phi^{n+1}(g))(c_0) = g(c_{-n-1})$ divides $h_1(c_0)$.  Each $r_i \in 
(T_B)_{\geq 1}$ and so $r_i(c_0) = 0$, and by assumption $\theta_{n+1} \not \in T_B$ 
and so $\theta_{n+1}(c_0) \in k^{\times}$.  Thus $g(c_{-n-1})$ divides 
$\Omega_{-n-1}$, which implies that $f(c_{-n-1})$ is a unit in $B$.  This contradicts property 
(2) above for $n \gg 0$.  Thus for $n \gg 0$ we must have 
$t_{n+1} \not \in \sum_{i = 1}^n t_i T_B$.  We 
conclude that $\sum t_i T_B$ is an infinitely generated 
right ideal of $T_B$, so $T\otimes_k B$ is not right 
noetherian and $T$ is not strongly right noetherian.  
\end{proof}

\section{Tensor Products of algebras}
In Proposition~\ref{strong noeth} we showed explicitly that $T$ is not strongly right 
noetherian by exhibiting a commutative noetherian $k$-algebra $B$ such that $T 
\otimes_k B$ is not right noetherian.  Necessarily, such a $B$ is not a finitely generated commutative algebra.  
By contrast, if we allow ourselves to tensor by noncommutative rings then we may find a   
finitely generated noetherian $k$-algebra $B'$ such that $T \otimes_k B'$ is 
not right noetherian.  In fact, we will see in the next theorem that one may take $B'$ to be $T$ itself.

In order to stay within the class of $\mb{N}$-graded algebras, 
in addition to tensor products it will be useful also to consider \emph{Segre products}, defined as follows. 
If $A$ and $B$ are two $\mb{N}$-graded algebras we let $A \segre_k B$ be the $\mb{N}$-graded algebra 
$\bigoplus_{n = 0}^{\infty} A_n \otimes_k B_n$.  The following lemma is then elementary.
\begin{lemma}
\label{segre vs tensor}
Let $A$ and $B$ be $\mb{N}$-graded algebras.  If $A \otimes_k B$ is left (right) noetherian, 
then $A \segre_k B$ is left (right) noetherian.
\end{lemma}
\begin{proof}
Since any homogeneous left ideal $I$ of $A \segre B$ satisfies $(A \otimes B)I \cap (A \segre B) 
= I$, a proper ascending chain of homogeneous left ideals of $A \segre B$ induces a 
proper ascending chain of left ideals of $A \otimes B$.
\end{proof}

We thank James Zhang for pointing out to us the following useful fact.  
\begin{lemma}
\label{finite pres} Let $A$ be connected $\mb{N}$-graded and noetherian.  Then $A$ is 
finitely presented.
\end{lemma}
\begin{proof}
Let 
\[
\dots \to \bigoplus_{i = 1}^{r_1} A[-d_{1i}] \to \bigoplus_{i = 1}^{r_0} A[-d_{0i}] \to A \to k \to 0
\]
be a graded free resolution of $_A k$ by free modules of finite rank.  Then one 
may check that $A$ has a presentation with $r_0$ generators and $r_1$ relations.
\end{proof}

The following theorem shows that it is possible to find two connected graded 
noetherian rings whose tensor product is noetherian on one side only, as well a pair 
of connected graded noetherian rings whose tensor product is noetherian on neither side.
\begin{theorem}
\label{tensor idealizer} 
Assume Hypothesis~\ref{second assump}, and in addition that $\{ c_n \}_{n \in \mb{Z}}$ is 
critically dense.  Let $T' = T \segre_k T^{op}$.  Then 
\begin{enumerate}
\item $T$ and $T'$ are noetherian finitely presented connected graded $k$-algebras.
\item $T \otimes_k T$ is left noetherian, but not right noetherian.
\item $T' \otimes_k T' \cong T' \otimes_k (T')^{op}$ is neither left nor right noetherian.
\end{enumerate}
\end{theorem}
\begin{proof}
(1) 
The ring $T$ is noetherian by Propositions~\ref{left noeth} and \ref{main noeth and chi result}.
In fact, by Proposition~\ref{strong noeth} $T$ is universally left 
noetherian. It follows immediately that $T^{op}$ is universally right noetherian.  
Thus $T \otimes_k T^{op}$ is both left and right noetherian.  By 
Lemma~\ref{segre vs tensor}, $T'$ is noetherian.  Then by 
Lemma~\ref{finite pres}, both $T$ and $T'$ are finitely presented.

(2) As we saw in part (1), $T$ is universally left noetherian, so that $T \otimes_k T$ is 
left noetherian.  Now we will prove that $T \otimes T$ is not right noetherian.  By 
Lemma~\ref{segre vs tensor}, it is enough to prove that $T \segre T$ 
is not right noetherian.  

For a graded ring $A$ we will use the abbreviation $A^s = A \segre A$.
Now let $X = \cproj U^s  \cong \mb{P}^d \times \mb{P}^d$.  The graded ring $U^s$ has 
the automorphism $\phi \otimes \phi$ with corresponding automorphism $\varphi \times 
\varphi$ of $X$.  The graded ring $S^s$ may be thought of as the left Zhang twist of 
$U^s$ by $\phi \otimes \phi$, and we identify the underlying vector spaces.  In 
particular, any homogeneous element of $S^s$ defines a vanishing locus in $X$.  Now 
let $\Delta \subset X$ be the diagonal subscheme, and let $J$ be the left ideal of 
$S^s$ consisting of those elements which vanish along $\Delta$.  Since $(\varphi 
\times \varphi)(\Delta) = \Delta$, it follows easily that $J$ is a two-sided ideal of 
$S^s$.  Writing $K = I \segre_k I$, a left ideal of $S^s$, we have $T^s = k \oplus K$.  Then to 
prove that $T^s$ is not right noetherian, by Lemma~\ref{right noeth} it will be enough 
to show that $(J \cap K)/JK$ is not finite-dimensional over $k$.   

Let $\circ$ indicate multiplication in the commutative ring $U^s$.  Since $J$ is 
invariant under $\phi \otimes \phi$, we have $J \circ K = JK$, and so it will 
be equivalent to prove that $M = (J \cap K)/(J \circ K)$ is not a torsion $U^s$-module.  To show 
this, we consider the corresponding sheaf $\wt{M}$ on $X$, look locally 
at the point $p = (c,c)$, and prove that $\wt{M}_p \neq 0$.  

Choose local affine coordinates
$u_1, \dots u_d$ for a principal open set $\mb{A}^d \subseteq \mb{P}^d$ such that the point $c$
corresponds to the origin.  Let $v_1, \dots v_d$ be the same coordinates for the equivalent
open set $\mb{A}^d$ in the second copy of $\mb{P}^d$, so that 
$u_1, \dots u_d, v_1, \dots v_d$ are local coordinates for an affine neighborhood 
$\mb{A}^{2d}$ of $p$ in $X$ such that $p$ is the origin in these coordinates.    
Now let $\mf{p}$ be the homogeneous 
prime ideal of $U^s$ corresponding to the point $p = (c, c)$. 
Setting $U' = (U^s)_{(\mf{p})} = \mc{O}_{X,p}$, $J' = J_{(\mf{p})}$, and 
$K' = K_{(\mf{p})}$, we have 
\[ 
\wt{M}_p = M_{(\mf{p})} \cong (J' \cap K')/ (J'K') 
\]
where we revert to the use of juxtaposition to indicate multiplication in the commutative local ring $U'$.
Explicitly, $U'$ is the polynomial ring $k[u_1, \dots u_d, v_1, \dots v_d]$ localized at the 
maximal ideal $\mf{m} = (u_1, \dots, u_d, v_1, \dots, v_d)$, $J' = (u_1- v_1, \dots, u_d 
-v_d)$, and $K' = (u_1, u_2, \dots, u_d)(v_1, v_2, \dots v_d)$.  Now it is clear that 
$w = u_1 v_2 - u_2 v_1 \in J' \cap K'$, but $w \not \in J'K'$ since $w \not \in \mf{m}^3 \supseteq J'K'$.
Thus $\wt{M}_p \neq 0$, as we needed to show.  

(3) Note that $(T')^{op} \cong T^{op} \segre T \cong T'$.  The fact that 
$T' \otimes T'$ is neither left nor right noetherian follows immediately from part (2).
\end{proof}

\begin{remark}
Assuming the setup of Hypothesis~\ref{second assump}, the ring $R = k \langle I_1 
\rangle \subseteq S$ which is generated by the degree $1$ piece of 
$T$ is a graded ring of the type studied in the article \cite{Ro1}.  In case the 
points $\{ c_n \}_{n \in \mb{Z}}$ are critically dense, this ring $R$ has similarly 
strange properties under tensor products.  For example, a similar but slightly more 
complicated version of the argument in Theorem~\ref{tensor idealizer}(2) above would 
show that $R \otimes_k R$ is neither left nor right noetherian.
\end{remark}

\section{proof of the main theorem}
In the final section, we recapitulate all of our preceding results to prove 
Theorem~\ref{main theorem}, which we restate as Theorem~\ref{main theorem redux} 
below.  The only thing we have left to show is that given the setup of 
Hypothesis~\ref{second assump}, there exists a plentiful supply 
of choices of a point $c \in \mb{P}^d$ and an automorphism $\varphi \in \aut 
\mb{P}^d$ such that $\mc{C} = \{ c_n \}_{n \in \mb{Z}}$ is critically dense.  
This situation has already been studied in the paper \cite{Ro1}; we repeat the 
result for the reader's reference as the next proposition.

We call a subset of a variety $X$ \emph{generic} if its 
complement is contained in a countable union of closed subvarieties $Z \subsetneq X$.  Note 
that as long as the base field $k$ is uncountable, any generic subset 
is intuitively ``almost all'' of $X$, in particular it is nonempty.  Thus the first part of the 
following proposition shows that if the base field $k$ is uncountable, then any 
suitably general pair $(\varphi, c)$ will lead to a critically dense set $\mc{C}$.  
The second part shows that in case $\cha k = 0$ we may easily write down many 
explicit examples of pairs $(\varphi, c)$ for which $\mc{C}$ is critically dense.    
\begin{proposition} \cite[Theorem 12.4, Example 12.8]{Ro1}
\label{crit dense result} Assume Hypothesis~\ref{second assump} and set $\mc{C} = 
\{c_n\}_{n \in \mb{Z}}$.  
\begin{enumerate}
\item Let $k$ be uncountable.  For any given $c \in \mb{P}^d$, there is a generic  
subset $Y \subseteq \aut \mb{P}^d = \pgl(k, d)$ such that if $\varphi \in Y$ then 
$\mc{C}$ is critically dense.
\item If $\cha k = 0$, $c = (1:1: \dots :1)$, and $\varphi$ is defined by 
\[
(a_0: a_1: \dots : a_d) \mapsto (a_0: p_1a_1: p_2a_2: \dots: p_d a_d), 
\] 
then 
$\mc{C}$ is critically dense if and only if $p_1, \dots, p_d$ generate a 
multiplicative subgroup of $k^{\times}$ which is isomorphic to $\mb{Z}^d$.  
\end{enumerate}
\end{proposition}

Finally, we summarize all of the properties that the ring $T$ has in case the set of points 
$\mc{C}$ is critically dense. 
\begin{theorem}
\label{main theorem redux} Assume Hypothesis~\ref{second assump}.  Let $k$ be 
uncountable and assume that the pair $(\varphi, c)$ is chosen so that $\mc{C} = \{ 
c_n \}_{n \in \mb{Z}}$ is critically dense.  Then the idealizer ring $T = \mb{I}(I) = 
T(\varphi, c)$ is a noetherian connected finitely presented graded ring with the 
following properties:
\begin{enumerate}
\item $T$ is left universally noetherian, but not strongly right noetherian.
\item $T \otimes_k T$ is 
left noetherian but not right noetherian.  The Segre product $T' = T \segre_k T^{op}$ is
also a finitely presented connected graded noetherian ring, but $T' \otimes_k T'$ is noetherian on neither side.
\item $\rProj T$ and $T\Proj$ have the same underlying category but non-isomorphic distinguished 
objects; specifically, $\rProj T \cong (\rQch \mb{P}^d, \mc{O}_{\mb{P}^d})$ and 
$T\Proj \cong (\rQch \mb{P}^d, \mc{I})$, where $\mc{I}$ is the sheaf of ideals corresponding 
to the point $c \in \mb{P}^d$. 
\item $T$ satisfies left $\chi_{d-1}$ but not left $\chi_d$, and $T$ fails $\chi_1$ on the right.
\item $\cd (\rProj T) = \cd (T\Proj) = d$.
\item Although no Veronese ring of $T$ is generated in degree $1$, one has isomorphisms 
$T\Proj \cong T^{(n)}\Proj$ and $\rProj T \cong \rProj T^{(n)}$ for all $n \geq 1$.
\end{enumerate}
\end{theorem}
\begin{proof}
Note that by Proposition~\ref{crit dense result}, we may indeed find a pair $(\varphi, c)$ 
so that $\mc{C}$ is critically dense.  Then $T$ is noetherian by 
Propositions~\ref{main noeth and chi result}(2) and \ref{left noeth}, and $T$ is 
finitely presented by Lemma~\ref{finite pres}.

Now (1) follows from Proposition~\ref{strong noeth}, and (2) from Theorem~\ref{tensor idealizer}.  

For (3), note that since $S$ is a left Zhang twist of $U$, we have $S\Gr \simeq U\Gr$ and so it easily 
follows that $S\Proj \cong U\Proj$.  Now the opposite 
ring $S^{op}$ of $S$ is isomorphic to the left Zhang twist of $U$ by $\phi^{-1}$; this may be 
checked directly, or see the proof of \cite[Lemma 4.2(1)]{Ro1}.  Thus we also have an isomorphism
$\rProj S \cong \rProj U$.  By Serre's theorem, we also have an equivalence of categories
$U\Qgr \simeq \rQch \mb{P}^d$, where $\rQch \mb{P}^d$ is the category of quasi-coherent sheaves on $\mb{P}^d$.  

Now using Proposition~\ref{cat equivs}, it follows that 
\[
T\Proj \cong (S\Qgr, \pi I) \cong (\rQch \mb{P}^d, \mc{I})
\]
and 
\[
\rProj T \cong (\rQgr S, \pi S) \cong (\rQch \mb{P}^d, \mc{O}_{\mb{P}^d}).
\]
Since $d \geq 2$, the ideal sheaf $\mc{I}$ which defines the closed point $c$ is not 
locally free, so in particular we have $\mc{I} \not \cong \mc{O}_{\mb{P}^d}$ and (3) 
is proved.  

Next, result (4) is 
a combination of Propositions~\ref{right chi} and ~\ref{main noeth and chi result}(1).
Since $S\Proj \cong U\Proj$ and $\rProj S \cong \rProj U$, it follows easily 
that $\cd(S\Proj) = \gd(S\Qgr) = \cd(\rProj S) = \gd(\rQgr S) = d$, and 
so (5) is a consequence of Proposition~\ref{cd}.  Finally, (6) follows from 
Proposition~\ref{Veroneses} and Lemma~\ref{basic props}(2).
\end{proof}

We close with a few remarks concerning Theorem~\ref{main theorem redux}. 
\begin{remark}
Theorem~\ref{main theorem redux}(2) shows that the tensor product of two noetherian finitely 
presented connected graded algebras (over an algebraically closed field) can fail to be noetherian.  This 
answers \cite[Appendix, Open Question $16'$]{GW}.  
\end{remark}

\begin{remark}
Suppose that $A$ is a connected graded noetherian ring satisfying left $\chi_1$ such that 
$A\Proj \cong (\rQch X, \mc{O}_X)$ for some proper scheme $X$.  
Keeler showed that in this case $A$ must be equal in large degree 
to a twisted homogeneous coordinate ring $B(X, \mc{L}, \sigma)$ where $\mc{L}$ is $\sigma$-ample 
\cite[Theorem 7.17]{Ke2}.  In particular, $A$ must be universally noetherian and must satisfy $\chi$ on both sides.  

Now consider instead connected graded noetherian rings 
$A$ with left $\chi_1$ such that $A\Proj \cong (\rQch X, \mc{F})$ 
for some proper scheme $X$, but where $\mc{F}$ is not assumed to be 
the structure sheaf.  Then $A = T$, where $T$ satisfies the conclusions of Theorem~\ref{main theorem redux}, 
is an example showing that rings with much more unusual behavior
may occur in this case.
\end{remark}

\section*{Acknowledgments}
Much of the research for this article was completed at 
the University of Michigan and the University of Washington; the author thanks 
both institutions.  The author also thanks James Zhang and Paul Smith for helpful conversations, and 
especially Toby Stafford for his many useful suggestions which improved this article.

\def\cprime{$'$}
\providecommand{\bysame}{\leavevmode\hbox to3em{\hrulefill}\thinspace}
\providecommand{\MR}{\relax\ifhmode\unskip\space\fi MR }
\providecommand{\MRhref}[2]{%
  \href{http://www.ams.org/mathscinet-getitem?mr=#1}{#2}
}
\providecommand{\href}[2]{#2}

\end{document}